\documentclass[12pt]{amsart}

\newcounter{defcounter}
\usepackage{enumitem,latexsym,amsthm}
\usepackage{amssymb,latexsym,eufrak,amsmath,amscd,graphicx}
\usepackage[all]{xy}
\usepackage{comment}
\usepackage{pdfsync}

\numberwithin{equation}{section}

\theoremstyle{plain}
\newtheorem{theorem}{Theorem}
\newtheorem{proposition}[theorem]{Proposition}
\newtheorem{corollary}[theorem]{Corollary}
\newtheorem{lemma}[theorem]{Lemma}

\newtheorem{proposition.definition}[theorem]{Proposition/Definition}

\newtheorem{theoremalpha}{Theorem}

\newtheorem{propositionalpha}[theoremalpha]{Proposition}

\theoremstyle{definition}

\newtheorem{remark}[theorem]{Remark}

\newcommand{\lra}{\longrightarrow}

\newcommand{\noi}{\noindent}
\newcommand{\PP}{\mathbf{P}}

\newcommand{\ZZ}{\mathbf{Z}}
\newcommand{\CC}{\mathbf{C}}

\newcommand{\OO}{\mathcal{O}}

\newcommand{\rk}{\text{rank }}

\newcommand{\HH}[3]{H^{{#1}} \big( {#2} , {#3} \big)}

\newcommand{\hh}[3]{h^{{#1}} \big( {#2} , {#3} \big)}

\newcommand{\AAA}{\mathbf{A}}

\newcommand{\ord}{\textnormal{ord}}

\newcommand{\pr}{\prime}

\newcommand{\Div}{\text{Div}}
\newcommand{\num}{\equiv_{\text{num}}}
\newcommand{\lin}{\equiv_{\text{lin}}}

\newcommand{\dra}{\dashrightarrow}

\newcommand{\Jac}{\textnormal{Jac}}

\newcommand{\divv}{\textnormal{div}}

\newcommand{\Linser}[1]{| \mspace{1.5mu} {#1} \mspace{1.5mu} |}
\newcommand{\linser}[1]{\Linser{{#1}}}

\newcommand{\Sym}{\textnormal{Sym}}

\newcommand{\ol}[1]{\overline{#1}}

\newcommand{\Tot}{\operatorname{Tot}}

\newcommand{\mmu}{\boldsymbol{\mu}}
\numberwithin{theorem}{section}

\newcounter{parasection}[section]
\renewcommand{\theparasection}{\thesection.\arabic{parasection}}

\newcommand{\parasection}[1]{%
  \refstepcounter{parasection}%
  \medskip
  \noindent{\itshape \theparasection. #1.}\quad
}

\begin{document}

\title{Some normal abelian surfaces in $\PP^3$}

\author{Nathan Chen}
\address{Instituto de Matemática Pura e Aplicada (IMPA)}
\email{nathan.chen@impa.br}

\author{Robert Lazarsfeld}
\address{Stony Brook University and University of Pennsylvania}
\email{robert.lazarsfeld@stonybrook.edu}

\author{Federico Moretti}
\address{Stony Brook University}
\email{federico.moretti@stonybrook.edu}

\begin{abstract} Every complex projective surface is birational to a singular surface $S \subseteq \mathbf{P}^3$, but it is unclear whether one can ask that $S$ has only isolated singularities. In particular, it was a folk question as to whether an abelian surface could admit such a model. We present examples to show that there do exist normal abelian surfaces in $\mathbf{P}^3$.  One of these examples deforms in a one-parameter family. \end{abstract}

\maketitle

\setlength{\parskip}{.13in minus .03in}

\section*{Introduction}

\renewcommand{\theequation}{\arabic{equation}}
\setcounter{equation}{0}

The purpose of this note is to exhibit some  abelian surfaces that are birational to surfaces in  $\PP^3$ having only isolated singularities.

Recall that any smooth complex projective  surface $X$ is birational to a surface $S \subseteq \PP^3$, but typically $S$ will be singular along a curve. Over the years there has been some interest in looking for realizations where $S$ has only isolated singularities, or in other words where $S$  is normal. Examples with non-zero irregularity $q \ne 0$ abound \cite{Zariski, Libgober, Oka},\footnote{For instance Zariski's famous example of a six-fold cyclic cover of $\PP^2$ branched along a sextic with six cusps lying on a conic, whose smooth model has $q = 1$, gives the surface $W^6 = Q_2(X, Y, Z)^3 - F_3(X, Y, Z)^2$, where $Q, F$ have degrees two and three respectively. } and   Hacking observed that work of Tokunaga \cite{TokunagaKodai, TokunagaPAMS}  yields  an example having Albanese dimension $= 2$. However it was a folk question (cf.\ \cite[Conjecture on p.\ 369]{NakamuraUmezu}) whether  abelian surfaces could  arise in this manner: the main result of \cite{NakamuraUmezu} is that  in any event there are no examples of degree $\le 5$. We show here that there do exist normal abelian surfaces in $\PP^3$.

Specifically, we prove
\begin{propositionalpha} \label{Intro.Elliptic.Curves.Example}
Let $E$  be the elliptic curve $   y^2 \, = \,  x^3 + 1 $ with $j = 0$. Then $E \times E$ is birational to a surface $S_9 \subseteq \PP^3$ of degree $9$ with only isolated singularities. In fact, $S_9$ can be given by the equation
\[
W^3 Q^3  \ = \ (X^3 - Y^3)(Y^3 - Z^3)(Z^3 - X^3)
\]
where $Q = Q(X, Y, Z)$ is a generic form of degree $2$.
\end{propositionalpha}
\noi The self-product of the elliptic curve $y^2 = x^3 - x$ also admits a normal model in $\PP^3$. 

\begin{theoremalpha} \label{Jacobian.Thm.Intro}
Let $C$ be the genus two curve 
\[ y^2 \ = \ x^5 + 1, \] and let $A = \Jac(C)$ be the Jacobian of $C$. Then $A$ is birational to a normal surface $S_{10} \subseteq \PP^3$ of degree $10$.  
\end{theoremalpha}
\noi In fact, we will see that there is a one-parameter family $\{ C_t \}$ of deformations of $C$ whose Jacobians are also birational to normal decics in $\PP^3$. 

 Unfortunately  we have nothing to say about the  arguably more interesting problem of finding obstructions to realizing a given smooth surface birationally in $\PP^3$  with only isolated singularities. In fact, to the best of our knowledge there are no surfaces  known that provably do not admit such models.

The starting point for these results is an elementary observation about  surfaces birational to an $n$-fold cyclic covering
\[
 f: V  \lra \PP^2. 
\]
Namely, assume that $f$ is branched along a \textit{reduced} plane curve of degree $nb$  in the sense  that $V$ is given by a  weighted homogeneous equation of the form
\begin{equation} \label{Cyclic.Cover.Equation}
T^n \ = \ B_{nb}(X, Y, Z)
\end{equation}
with $B$ square-free. Then  $V$ is birational to a normal surface in $\PP^3$ (Proposition \ref{Constructing.Normal.Surfaces}). For the product of the elliptic curves appearing in Proposition \ref{Intro.Elliptic.Curves.Example} such a representation essentially appears in the literature \cite{RU}, and the result follows without further ado. (We give a quick   derivation as warm-up for the constructions in \S \ref{Genus.2.Jacobian.Section}, and treat also the elliptic curve with $j = 1728$.)

However many of the classical  cyclic multiple planes do not have this form over the line at infinity,  so this picture does not apply directly. This is in particular the case in the interesting work of  Artal Bartolo, Cogolludo-Agust\'in, and Libgober \cite{ACL}. These authors construct a ten-sheeted cyclic covering $A_1 \lra \PP^2$ branched over a carefully chosen curve with an $A_4$ singularity, and then identify $A_1$ as birational to 
 the abelian surface $A = \Jac(C)$. But as they observe, the covering in their example factors as a composition
 \[
 A_1 \underset{\text{deg} \, 5} \lra R \underset{\text deg \, 2}\lra \PP^2.
 \]
 of two coverings of degrees $5$ and $2$.
One sees that $R$ is rational, and AI produced a birational identification of $R$ with $\PP^2$ that brought the first covering into the form \eqref{Cyclic.Cover.Equation} with $n = 5$ and $b = 2$. This in principle  proves Theorem \ref{Jacobian.Thm.Intro}, but the specific form of the presentation in \cite{ACL} is not very convenient for our purposes. Instead, we work here from the top down. Specifically, we show   that the quotient map for the evident $\mu_5$ action on $A$ can be described birationally by a cyclic cover
\[
A_2 \lra \PP^1 \times \PP^1
\]
of degree $5$ branched along the union $B = \sum B_i$ of five curves $B_i \in \linser{\OO_{\PP^1 \times \PP^1}(1,1)}$. Projecting to $\PP^2$ yields the required representation.

 Many of the ingredients in the story described here were no doubt known in
one form or another, but we hope there is some interest in putting them
together in the present context. The connection between cyclic coverings and normal surfaces in projective
space seems to be  implicit in earlier work, although we are unaware of
any reference where Proposition \ref{Constructing.Normal.Surfaces} is stated
in this form. On the other hand, there is a large and beautiful literature
concerning the geometry of cyclic multiple planes. Besides the papers
already mentioned, we point to \cite{Zariski}, \cite{Libgober}, 
\cite{LibgoberQA}, \cite{Naie}, and
\cite{CataneseCiliberto}.

Concerning the organization of this note, we start in \S \ref{Cyclic.Covering.Lemma.Section} with a review of some general facts about cyclic coverings, and then explain how they can give rise to normal surfaces in $\PP^3$. As warm-up, we treat the self-product of elliptic curves with $j = 0$ and $1728$. Section \ref{Genus.2.Jacobian.Section}, where we study the Jacobian of the curve $C$ in Theorem \ref{Jacobian.Thm.Intro},  contains the main geometric work. Since the picture strikes us as being rather pretty, our exposition here is more leisurely than what is strictly required to establish the Theorem. Finally, in the third section we sketch very briefly how one can deform the construction while maintaining the normal model in $\PP^3$.

\subsection*{Acknowledgements and AI Use} We are grateful to Fabrizio Catanese, Paul Hacking, J\'anos Koll\'ar and Rita Pardini for valuable communications and correspondence. We have made essential use of ChatGPT   at every stage of this project; in particular it suggested the constructions explained in Sections \ref{Genus.2.Jacobian.Section}  and \ref{Deformation.Section}.  We also consulted Claude on a couple of occasions.  The authors of course take responsibility for the presentation and accuracy of the material.


\section{Cyclic coverings and normal surfaces} 
\label{Cyclic.Covering.Lemma.Section}

\renewcommand{\theequation}{\thesection.\arabic{equation}}

\subsection{Generalities on cyclic coverings} \label{Cyclic.Generalities}

We start with some remarks  about cyclic coverings that will clarify later constructions. The picture itself is well-known (cf \cite{Catanese2012}, \cite{Pardini}), but it will be convenient to set down the facts from the perspective we need.  For concreteness we
work with surfaces, but everything we say goes over to higher dimensions with essentially no  modifications. 

By an \(n\)-fold cyclic cover we understand a finite mapping
\begin{equation} \label{Cyclic.Cover.Setup}
f : V \lra U
\end{equation}
of (possibly non-complete) surfaces such that the extension of function fields
\[
\CC(V) \, \supset \, \CC(U)
\]
is cyclic of degree \(n\). Given $f$, we   fix once and for all a generator
\[
\sigma\ \in \ \operatorname{Gal}\big(\CC(V)/\CC(U)\big).\]
  We will always assume that \(V\) is normal and that \(U\) is at least
normal, and usually smooth. We will also generally  consider such coverings only
up to birational equivalence. Thus  if \(U_1\sim_{\mathrm{bir}} U\), we get
a corresponding model
\[
f_1:V_1\lra U_1
\]
of the covering by making the identification
\[
\CC(U)=\CC(U_1)
\]
and taking \(V_1\) to be the normalization of \(U_1\) in \(\CC(V)\). Therefore  in practice we are free to replace \(U\) by a desingularization.

An \(n\)-fold cyclic cover \eqref{Cyclic.Cover.Setup} determines---and in
happy cases is determined by---branching data along the prime divisors
of \(U\). We will call these the \textit{Kummer multiplicities}  of the cover, and they are
defined as follows.\footnote{There does not seem to be a standard name for the quantities in question, although they appear in many computations with  cyclic coverings. Catanese refers to them as local monodromies. They carry  the same information as the residues of the logarithmic connections studied by Esnault and Viehweg \cite{EV}, and in the present cyclic setting they are equivalent to knowing the building data of the covering introduced by Pardini \cite{Pardini}.}  Since we work over \(\CC\), Kummer theory applies. Set $\zeta = \zeta_n = \exp(\frac{2\pi i}{n})$. Then there exists  an element \(\phi  \in \CC(V)\), with $
\sigma(\phi)\, = \, \zeta \phi
$ and 
\[ \phi^n \ = \ b \, \in \, \CC(U)^*, \] such that
$
\CC(V) = \CC(U)(\phi).$
 The element \(\phi\) is determined up to
multiplication by an element in  \(\CC(U)^*\).

Now let \(E\subseteq U\) be a prime divisor. We define the Kummer multiplicity of the
cover along \(E\) to be   \[
k_E(f) \, 
=
\, \operatorname{ord}_E(b)\]
considered as an integer (mod $n$). 
This is a well-defined element in $\ZZ/n\ZZ$:   if we replace \(\phi\) by
\(g\phi\), where \(g\in \CC(U)^*\), then \(b\) is replaced by \(g^n b\), and
therefore
\[
\operatorname{ord}_E(g^n b)
=
n\operatorname{ord}_E(g)+\operatorname{ord}_E(b).
\]
The cover is unramified at the generic point of \(E\) if and only if
$k_E(f) = 0$. More generally, the ramification index along \(E\) is given by 
\[
\frac{n}{\gcd(n, k_E(f))}.
\]
In particular, when \(n\) is prime, every nonzero multiplicity indicates total
ramification along \(E\).

It is convenient to assemble these multiplicities into a divisor on $U$. Given a cyclic covering \eqref{Cyclic.Cover.Setup}, define
\[
B \, = \, B_f \, = \, \sum_E \, k_E(f) \cdot E \, \in \,\Div(U) _{\ZZ/n\ZZ}.
\]
We call this the (Kummer) \textit{branch divisor} of the covering. Thus $B$ has $\ZZ/n\ZZ$ coefficients, but in concrete calculations we represent it as an integer divisor with the understanding that the coefficients are only well-defined (mod $n$). Observe also that this is not the usual branch divisor of a covering $f : V \lra U$,  ie the push-forward of ramification divisor of $f$: rather the picture is that if one thinks of defining a cyclic covering by an equation $t^n = b(x,y)$, then $B = \divv(b)$. 

These divisors are useful for keeping track of the Kummer multiplicities under  changes of birational models:
\begin{lemma} \label{Branch.Transformation.Rule}
Let $\mu : U \lra U_1$ be a proper birational morphism, with $U$ and $U_1$ non-singular, and let
\[
f : V \lra U \ \ , \ \ f_1 : V_1 \lra U_1
\]  be the corresponding models of a cyclic covering. Denote by $B$ and $B_1$ the Kummer branch divisors of $f$ and $f_1$ respectively. Then 
\[
B_1 \ = \ \mu_* B \ \ , \ \ B \ = \ \mu^* B_1.
\]
\end{lemma}
\begin{proof} This follows from the facts that the field extension has not changed and that $\ord_E$ is a valuation on $\CC(U) = \CC(U_1)$ for  any divisor $E$ lying over $U$ or $U_1$.
\end{proof}

We will say that a cyclic covering $f : V \lra U$ has \textit{reduced branch divisor} if all the non-zero coefficients of $B_f$ are $\equiv 1 \pmod{n}$, so that it can be represented by a reduced integral divisor. This condition will allow us to write an explicit equation defining $V$, but first we require a little notation. 

Assume that $U$ is non-singular,  let $L$ be a line bundle on $U$, and denote by $\pi : \mathbf{L} \lra U$  its total space. Recall that there is a tautological section \[ T \  \in \ \HH{0}{\mathbf{L}}{\pi^* L}\] that serves as a coordinate in each fibre of $\pi$ (cf \cite[Proof of Proposition 4.1.6]{PAG}). Hence if $s \in \HH{0}{U}{L^n}$ is any section, then 
$
T^n - \pi^*(s) $ is a section of $\pi^* L^n$.
  To lighten the discussion, we call its divisor  the hypersurface in $\mathbf{L}$ defined by $T^n = s$.
\begin{lemma} \label{Reduced.Covering.Lemma}
Let \[       f : V \lra U
\] be a cyclic covering as above, with $U$ non-singular.
Suppose that the Kummer
branch divisor of \(f\) is reduced. Then there exist a line bundle \(L\) on
\(U\) and  a section
\[
        s \ \in\ H^0(U,L^n)\]
such that $\divv(s)$ is a reduced integral representative of $B_f $, and moreover \(V\) is
isomorphic over \(U\) to the hypersurface in $\mathbf{L}$ given by
\begin{equation} \label{Cyclic.Cover.Equation.in.L}
        T^n \ =\ s.
\end{equation}
\end{lemma}

\begin{remark}  Absent the reducedness hypothesis, one can still    find $ s  \in H^0(U,L^n)$ whose divisor represents $B_f$, but  then one generally has to normalize the hypersurface in $\mathbf{L}$ defined by  \eqref{Cyclic.Cover.Equation.in.L} to get $V$. 
\end{remark}

\begin{proof}
Choose a Kummer generator $ \phi \in \CC(V)$ with 
$    \phi^n = b   \in \CC(U)^*$, and  let  \(B\) be a reduced effective  integer divisor representing  the Kummer branch divisor. Then by definition
\[
        \operatorname{div}(b) \ =\ B + nD
\]
for some divisor \(D\) on \(U\). Thus
$B \lin -nD$. Set
\[
        L \ =\ \OO_U(-D).
\]
Then $        L^n \cong  \mathcal O_U(B)$, and we take
  \[     s \ \in\ H^0(U,\mathcal O_U(B)) \ \cong \ H^0(U,L^n)
\]
to be the section whose zero-divisor is \(B\).\footnote{More precisely, we take $s=b$, viewed as a section of
$L^{\otimes n}=\OO_U(-nD)$. Its zero-divisor is $B$,
and this choice ensures that the equation $T^n=s$
defines the original Kummer extension.}
It remains only to show that $V$ is given by the equation
\eqref{Cyclic.Cover.Equation.in.L} and this is where we use the hypothesis that $B$ is reduced. In fact, \eqref{Cyclic.Cover.Equation.in.L}  defines  a finite cyclic cover of \(U\) sitting in the total space $\mathbf{L}$ that determines the
same function field extension as \(f\). So it suffices to show that  \eqref{Cyclic.Cover.Equation.in.L} cuts out a normal variety, or equivalently -- since $\mathbf{L}$ is smooth -- that the hypersurface $T^n - s = 0$ is nonsingular in codimension one. 
But this is clear: since $s$ is reduced its singularities have codimension $\ge 2$ in $U$, and hence the singularities of $\{ T^n - s = 0\}$ have codimension $\ge 2$ in $\mathbf{L}$. \end{proof}

\begin{corollary} \label{Reduced.Covering.P2}
Let $f : V \lra \PP^2$ be a cyclic covering with reduced branched divisor $B$. Then $\deg B = bn$ for some $b \ge 1$ and $V$ is defined by a weighted homogeneous equation
\begin{equation} \label{Cyclic.Cover.P2}
T^n \ = \ F_{nb} \  \ \text{in } \ \PP(1,1,1,b),
\end{equation}
where $F_{nb} = F_{nb}(X,Y,Z)$ is a square-free homogeneous polynomial of degree $nb$. \qed \end{corollary}

\begin{remark} \label{Cyclic.Covering.Data}
In the situation of Lemma \ref{Reduced.Covering.Lemma}, the data needed to define a cyclic covering are the Kummer branch divisor $B$, an $n^\text{th}$ root $L$ of the line bundle $\OO_U(B)$, and a choice of isomorphism $L^{\otimes n}\cong \OO_U(B)$. The possible choices of $L$, when they exist, form a torsor under $\operatorname{Pic}(U)[n]$. For a fixed $L$, different choices of isomorphism can give different coverings, the ambiguity being a unit modulo $n^\text{th}$ powers. When $U$ is complete (so that every unit is an $n^\text{th}$ power) and $\operatorname{Pic}(U)$ has no $n$-torsion, the Kummer branch divisor determines the cyclic covering. This is the cyclic case of Pardini's building-data description of abelian covers \cite{Pardini}. In the language of this paper, having a reduced branch divisor is equivalent to being a simple cyclic covering \cite[Definitiion 1.3]{Pardini}. See also \cite[\S 1, 4a]{CataneseDiffDef}.\end{remark}

\subsection{Construction of normal surfaces} 
We now explain how cyclic coverings of $\PP^2$ with reduced branch divisors can be realized birationally as surfaces in $\PP^3$ having only isolated singularities.

\begin{proposition} \label{Constructing.Normal.Surfaces} Let
\[ f : V \lra \PP^2 \]
be an $n$-fold cyclic covering branched along a reduced curve $B \subseteq \PP ^2$ of degree $nb$ as in Corollary \ref{Reduced.Covering.P2}. Then $V$ is birational to a normal surface
\[ S \ = \ S_{bn} \ \subseteq \ \PP^3 \]
of degree $bn$.
\end{proposition}
\noi The proof of the Proposition is quick but perhaps not very intuitive. To compensate, we first describe briefly the geometry of a slightly different birational embedding of $V$.  We assume henceforth that $b \ge 2$, since if $b = 1$ then \eqref{Cyclic.Cover.P2} already exhibits $V$ as embedded in $\PP^3$.

Write $\OO_V(k) = f^* \OO_{\PP^2}(k)$ and note that $f^*$ determines an isomorphism
\[ \HH{0}{\PP^2}{\OO_{\PP^2}(k)}  \overset{\sim} \lra  \HH{0}{V}{\OO_V(k)} \ \ \text{when } k \le b-1,\]
whereas $\HH{0}{V}{\OO_V(b)} $ is spanned by $\HH{0}{\PP^2}{\OO_{\PP^2}(b)}$ together with (the restriction to $V$ of) the weight $b$ coordinate $T$ from Corollary \ref{Reduced.Covering.P2}. By a slight abuse of notation, we will think of  weighted homogeneous polynomials in $X, Y, Z$ and $T$ as defining sections of the appropriate line bundle on $V$.
Choose now general sections
\[
Q  = Q_{b-1}(X,Y,Z) \in \HH{0}{V}{\OO_V(b-1)}  \ \ , \ \ G = G_b(X,Y,Z,T) \in \HH{0}{V}{\OO_V(b)},
\]
defining curves $C, D \subseteq V$, and consider the four-dimensional subspace   $\Lambda \subseteq \HH{0}{V}{\OO_V(b)}$  spanned by the sections
\[ G \ , \  XQ \ , \ YQ \ , \ ZQ.
\] The rational map
\[ \varphi_{\linser{\Lambda}}:  V \dra \PP^3
\]
is birational onto its image, because $T$ appears linearly in $G$, and away from  $Q = 0$ we can recover $X, Y, Z$ and $G$ from the coordinates of  $\PP^3$.  

The base locus of $\linser{\Lambda}$ consists of the $b(b-1)n$ intersection points of $C$ and $D$. Blowing these up gives a morphism $\varphi^\pr : V^\pr \lra \PP^3$ whose image is a surface $S \subseteq \PP^3$ of degree $bn$. This morphism contracts the proper transform $C^\pr \subseteq V^\pr$ to a highly singular point of $S$, but we claim that $S$ is actually non-singular in codimension one. To see this it is enough to show that a general hyperplane section of $S$ is a non-singular plane curve, and this follows from a computation with the adjunction formula. Namely, it turns out  that the curves defining $\varphi^\pr$ -- which in any event give smooth models of this hyperplane section --   have the same genus as a smooth plane curve of degree $bn$.

An alternative  argument proceeds by simply writing down an equation for a slightly different model  of $V$ in $\PP^3$. 
\begin{proof} [Proof of Proposition \ref{Constructing.Normal.Surfaces}]
We keep the notation from the preceding discussion. Thus \[Q\ = \ Q_{b-1}(X,Y,Z)\] is a
general form of degree \(b-1\). This time we consider the slightly
different rational map
$
        V \dashrightarrow \PP^3
$
defined by the four sections
\[
        T,\quad XQ,\quad YQ,\quad ZQ
\]
of \(\OO_V(b)\). As before, this map is birational onto its image.
The plan is to write down the equation of the image and check directly that it is
non-singular in codimension one and hence normal. 

Let
$ [w,x,y,z]$ be homogeneous coordinates on the target \(\PP^3\). Since
\[
        [x,y,z] \ =\ [XQ,YQ,ZQ] \ =\ [X,Y,Z]
\]
away from \(Q=0\), and since
$  w \ =\ T, $
the defining equation
$T^n \ =\ F_{nb}(X,Y,Z)$
homogenizes  in the target coordinates to
$ w^n Q(x,y,z)^n \ =\ F_{nb}(x,y,z).$
Thus the image is the surface
\[
        S \ =\
       \Big\{
        w^n Q(x,y,z)^n \ =\ F_{nb}(x,y,z)
        \Big\}
        \ \subseteq\ \PP^3_{[w,x,y,z]} .
\]
In particular \(S\) has degree \(nb\), and it is birational to \(V\).

It remains to show that \(S\) is non-singular in codimension one. For
simplicity of notation, write $F = F_{nb}$ and set
\[
        \Phi \ =\ w^n Q^n-F
\]
Note first that $   {\partial \Phi / 
       \partial w}
        =
        nw^{n-1}Q^n,$
so \(S\) is non-singular away from \(wQ=0\). Suppose next that \(w=0\). Then
the equation of \(S\) gives
$ F (x,y,z) = 0,$ and the remaining partial derivatives restrict to the partial derivatives of
\(-F\). Thus any singular point with \(w=0\) lies over the finite singular
locus of the reduced plane curve
$       F  = 0.$
 Next suppose that \(Q=0\). On this locus, the equation gives
$F= 0$, and the partial derivatives of $\Phi$ with respect to \(x,y,z\) again reduce to those of $-F$. 
Hence a singular point with \(Q=0\) would have to lie over
\(\operatorname{Sing}(F=0)\). For general \(Q\), this does not occur,
except possibly at the point
$
       [1,0,0,0]
$
which is at worst an isolated singularity. Therefore the singular locus of
\(S\) is finite.  
\end{proof}

\begin{remark} \label{Monic.Hypersurface}
More generally, consider a surface
\[ V \, \subseteq \, \PP(1,1,1,b) \]
defined by a monic weighted homogeneous polynomial
\[
T^n \, + \, a_1T^{n-1} \, + \, \ldots \, + \, a_{n-1}T \, + \, a_n  \ = \ 0, 
\]
where $a_i = a_i(X,Y,Z)$ is homogeneous of degree $ib$. (We will call $V$ a ``monic hypersurface.") If $V$ is normal, then the same transformation as before yields a model $S \subseteq \PP^3$ of $V$ as a  surface in $\PP^3$ with only isolated singularities. This is established by either of the two arguments just presented; we leave details to the reader. 
\end{remark}

\begin{remark} \label{Higher.Dim.Extension.Prop.1.6}
Proposition \ref{Constructing.Normal.Surfaces} also extends to  all
dimensions. Suppose that
\[
        T^n \ =\ F_{nb}(X_0,\ldots,X_r)
\]
defines an \(n\)-fold cyclic cover of \(\PP^r\), where \(F = F_{nb}=0\) is
reduced. Choose a general form
$ Q=Q_{b-1}(X_0,\ldots,X_r),$
and consider the rational map given by
\[
        T,\quad X_0Q,\quad \ldots,\quad X_rQ .
\]
With homogeneous coordinates
$ [w,x_0,\ldots,x_r]$
on the target \(\PP^{r+1}\), the image is the hypersurface
\[
        w^n Q(x_0,\ldots,x_r)^n  \ = \ F (x_0,\ldots,x_r).
\]
If $     \operatorname{Sing}(F=0) $ has dimension $s$ 
then for \(Q\) general, the singular locus of this hypersurface also  has
dimension $s$
apart from the possible additional isolated point
$  [1,0,\cdots,0].$
In particular, since \(F=0\) is reduced, the resulting hypersurface is
non-singular in codimension one. The same conclusion also holds for any  normal monic hypersurface in $\PP(1^{r+1}, b)$ provided that $r \ge 2$.
\end{remark}


\subsection{First application: self-product of elliptic curves}  As warm-up, we give a quick proof of Proposition \ref{Intro.Elliptic.Curves.Example} and its companion for the elliptic curve $F = \{ y^2 = x^3 - x \}$, although we won't try to derive the specific equation stated in the Introduction. 

In view of  Proposition \ref{Constructing.Normal.Surfaces}, the  issue is to realize $E \times E$ and $F \times F$   birationally as cyclic covers of $\PP^2$ with reduced branch divisors. 
In fact, the following Lemma shows that even less is required.
\begin{lemma} \label{Rational.Quotient.Lemma}
Let $f_0 : V_0 \lra Y_0$ be a cyclic covering of  projective surfaces, with $V_0$ normal and $Y_0$ smooth. Assume that $f_0$ has reduced branch divisor, and that $Y_0$ is rational. Then $f_0$ is represented by a model
\[  f : V \lra \PP^2 \]
with reduced branch locus. 
\end{lemma}

\begin{proof}
Blow down ($-1$)-curves on $Y_0$ to construct a birational morphism
\[ \mu_0: Y_0 \lra Y_1 \]
where $Y_1$ is minimal. By Lemma \ref{Branch.Transformation.Rule}, the resulting model $f_1 : V_1 \lra Y_1$ of the covering still has reduced branch divisor. Being a minimal rational surface, $Y_1$ is either isomorphic to $\PP^2$ or to a Hirzebruch surface $\mathbf{F}_n$ for some $n \ge 0$. In the latter case, a sequence
of elementary transformations, followed by the contraction
of a $(-1)$-section, gives a birational map
$\mathbf{F}_n\dashrightarrow\PP^2$. By choosing the centers of the
blowings-up away from the branch divisor, we ensure that
all Kummer multiplicities remain either $0$ or $1$.
Then the resulting covering of $\PP^2$ has reduced branch
divisor, as required.  
\end{proof}

 \subsubsection*{Self-product of $E$} We now apply this to the self-product of the elliptic curve  
\[E \ : \, \  y^2=x^3+1 \] with \(j=0\). We will think of the point at infinity
\(\infty\in E\) as the origin on \(E\).   Setting \(\omega=e^{2\pi i/3}\),  there is a \(\mmu_3\)-action on
\(E\) given by
\[
        \omega\cdot (x,y)=(\omega x,y),
\]
which has three fixed points $p_+=(0,1)$, $p_- = (0,-1)\),  and \(\infty\). We now consider the diagonal
action on \(E\times E\): this is an order \(3\) automorphism with  \(9\) fixed points. 

Denote by 
\[
        \pi:E\times E\lra Y_0=_{\textnormal{def}} (E\times E)/\mmu_3.
\]
the quotient map. Each of the
\(9\) fixed points  gives rise to a cyclic quotient singularity of type
\(\frac{1}{3}(1,1)\) on $Y_0$;  this is local analytically isomorphic to a
cone over a rational normal cubic. The minimal resolution
\[
        \mu:Y\longrightarrow Y_0
\]
is achieved by blowing up each of the singular points \(y_i\in Y_0\) just
once, and over each \(y_i\) the exceptional divisor \(E_i\subset Y\) is
precisely a \((-3)\)-curve.
In coordinates, the \(\mmu_3\)-action in a neighborhood of a fixed point on
\(E\times E\) locally looks like
\[
        \omega\cdot (u,v)=(\omega u,\omega v)
\]
on \(\mathbf A^2_{u,v}\). One checks that in the resulting covering $f : V \lra Y$, the Kummer  multiplicity of each exceptional divisor \(E_i\) is  $k_{E_i} = 1$,  
so $f$ has reduced branch locus. 

It remains to show that $Y$ is rational, and this follows from Castelnuovo's criterion. In fact, since $Y_0$ has rational singularities, \[
        \HH{1}{Y}{\mathcal O_Y}
        =
        \HH{1}{E\times E}{\mathcal O_{E\times E}}^{\mmu_3}
        =
        0.
\]
Moreover, adjunction shows that each $E_i$ appears with coefficient $-1/3$ in $K_Y - \mu^* K_{Y_0}$, and hence $3K_Y \num- \sum E_i$. Therefore  all the plurigenera $Y$ vanish, and in particular $P_2(Y) = 0$, establishing rationality. Therefore Lemma \ref{Rational.Quotient.Lemma} applies and we see that $E \times E$ is realized birationally as a normal surface in $\PP^3$. 

\begin{remark} This quick argument ignores a lot of beautiful geometry explained in [14, \S 2, \S 3]. In brief, there are twelve distinguished copies of $E$ in $E\times E$ that are stable under the action of $\mmu_3$: start with the four curves
\[
T_0=\{0\}\times E,\quad T_1=\Delta_E,\quad
T_\infty=E\times\{0\},\quad T_{\zeta}=\Gamma_{\zeta},
\]
where $\Gamma_{\zeta}$ is the graph of the order-six automorphism $\zeta:E\to E $ given by $\zeta (x,y) = (\omega x, -y)$. Put $K=\{\infty,p_+,p_-\}\subseteq E$. Translating $T_0,T_1,T_{\zeta}$ by the points $(p,0)$, and $T_\infty$ by $(0,p)$, for $p\in K$, gives twelve copies of $E$, meeting only at the nine fixed points of the $\mmu_3$ action, with four passing through each fixed point. The images of these curves in $Y_0$ lift to  twelve disjoint $(-1)$-curves on $Y$, of which each exceptional divisor  $E_i$ meets four. Blowing these down yields a model $Y_1$ of $Y$ isomorphic to $\PP^2$, in which the images of the nine $E_i$ are the lines of the dual Hesse configuration. This leads to the equation
\[
W^3=(X^3-Y^3)(Y^3-Z^3)(Z^3-X^3)
\]
describing the covering. \end{remark} 

\subsubsection*{Self-product of $F$}  We now sketch how to apply a similar analysis to  the self-product of the elliptic curve
\[
F \ : \ \  y^2 = x^3 - x. \]
Consider the order four automorphism  of $F$  given by $\alpha(x,y) = (-x, iy)$, which has $\infty, 0 = (0,0) \in F$ as fixed points, and interchanges $p_{\pm} =_{\text{def}}  (\pm 1, 0)$ (which are fixed by $\alpha^2$). 
Write $\sigma = (\alpha, \alpha)$ for the diagonal action on $F \times F$, and let
\[ Y_0 \ = \ (F \times F) / \mmu_4\]
be the quotient. Then $Y_0$ has four cyclic quotient singularities of type $\tfrac{1}{4}(1,1)$ and six of type $\tfrac{1}{2}(1,1)$. Writing $\mu: Y \lra Y_0$ for the minimal resolution, we find four $(-4)$-curves $G_\alpha \subseteq Y$ and six $(-2)$-curves $J_\beta \subseteq Y$. The Kummer multiplicities are:
\[
k_{G_\alpha}\, = \,  1 \ \ , \ \ k_{J_\beta} \, = \, 2. 
\]
As before one sees that $Y$ is rational, but in order to apply Lemma \ref{Rational.Quotient.Lemma} we will need to get rid of the $J_\beta$.

To this end, consider the six  $\sigma$-stable copies of $F$ in $F\times F$:
\[
\{ \infty \} \times F \ , \ \{0\} \times F\ , \ F \times \{ \infty\} \ , \  F \times \{ 0 \} \ , \ \Delta \ , \ \Gamma_\alpha, 
\]
where $\Delta$ is the diagonal and $\Gamma_\alpha$ is the graph of $\alpha$. Their images on $Y$, which we denote by $\Gamma_i$ are pairwise disjoint $(-1)$-curves. Each $\Gamma_i$ meets two of the $G_\alpha$, and  one of the $J_\beta$. Blowing down the $\Gamma_i$ we arrive at a surface in which the image of each $J_\beta$ is a $(-1)$-curve. We may then blow these down to arrive at a model in which all the Kummer multiplicities $ = 1$.

 
\section{The genus two Jacobian} \label{Genus.2.Jacobian.Section}

This section is devoted to the construction underlying Theorem \ref{Jacobian.Thm.Intro}
 from the Introduction. Although quotients of abelian varieties by cyclic groups are classical objects
in the Bagnera--de Franchis tradition, cf. \cite{CataneseTopMethods}, we will
need here a rather explicit birational model of this particular quotient. 

In outline, there  is an evident $\boldsymbol{\mu}_5$-action on $A$, and the main result (Theorem \ref{Main.Thm.Quotient}) is that the quotient $A \lra A/\mmu_5$ has a canonical birational model 
\[
A^\pr \lra \PP^1 \times \PP^1
\]
with reduced branch divisor  $B = \sum B_i$ consisting of the union of five curves of type $(1,1)$.  Once we know this, Theorem \ref{Jacobian.Thm.Intro} follows immediately from the discussions in the previous section.  

As explained in \S \ref{Cyclic.Generalities}, all the important geometry takes place on the base: any birational model of the quotient determines a unique normal representative of the covering. So the issue is to understand how to pass from $A/\mmu_5$ to $\PP^1 \times \PP^1$ while keeping track of the Kummer branch multiplicities. There are a certain number of local assertions to check, and we have chosen to verify these by writing everything out explicitly. The main birational constructions start in Paragraph \ref{Birational.Parasection}, and the reader willing to invoke Lemma \ref{Rational.Quotient.Lemma} can skip most of Paragraph \ref{Conclusion.Parasection}.

\parasection{Set-up: $C$ and its Jacobian}  Our focus in this section is the genus two curve $C$ given in affine coordinates by
\begin{equation} \label{Affine.Eqn.C}
y^2 \ = \ x^5 \, + \, 1.
\end{equation} Near the point at infinity $\infty \in C$ we use the coordinates $x_\infty = 1/x$, $y_\infty = y/x^3$ so that the equation of $C$ becomes
\[ y_\infty^2 \ = \  x_\infty \, + \, x_\infty^6.\]
We fix the basis
\begin{equation} \label{Basis.of.Diffs}
\omega \, = \, \frac{dx}{y} \, = \, - x_\infty \frac{dx_\infty}{ y_\infty} \qquad \eta \, = \,  x \frac{dx}{y} \, = \, - \frac{dx_\infty }{ y_\infty}
\end{equation}
for the holomorphic differentials on $C$. Referring to the affine equation \eqref{Affine.Eqn.C}, we also single out two special points on $C$ that will play an important role later:
\[ p_+ \, = \, (0, 1) \ \ , \ \  p_- \, = \, (0, -1).
\]

Set $A = \Jac(C)$. Viewing the Jacobian as parameterizing divisor classes of degree zero, we take $\infty \in C$ as the base-point for defining the Abel-Jacobi embedding
\[ a : C \lra \Jac(C) \ \ , \ \ p \mapsto [ p - \infty]. \]
Thus $a(\infty) = 0$. We then fix global flat holomorphic coordinates $u , v$ on $A$ (or, more pedantically, on the universal covering of $A$) by requiring that 
\[
\omega \, = \, a^*(du) \ \ , \ \ \eta \, = \, a^*(dv). 
\]
We put $t_{\pm} = a(p_{\pm}) \in A$. Observe that these are five-torsion points on $A$ since $5(p_\pm - \infty) = \divv( y \mp 1)$, and $t_- = -t_+$.

\begin{figure} 
\centering
\includegraphics[scale = .75]{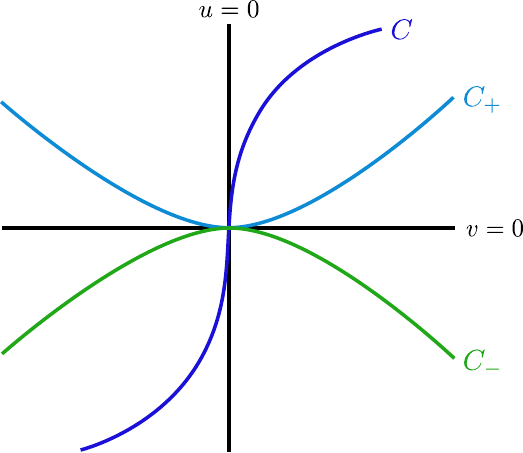}
\caption{Abel-Jacobi curves near 0.}
\label{AJCurves.Figure}
\end{figure}

We  will need to know something about the equations defining the Abel-Jacobi image  $C\subset A$ near  $0$ and $t_{\pm}$. 
\begin{lemma} \label{local.eqn.AJ.Curve}
In a neighborhood of $0$, $C \subseteq A$ is defined by a local analytic equation of the form
\[
u \ = \ c_1 \cdot v^3 \, + \, \textnormal{HOT}(v) \ \ \ (c_1 \ne 0).
\]
In  neighborhoods of $t_{\pm}$, the local analytic equations in coordinates centered at those points are \[
v \ = \ \pm c_2 \cdot u^2 \, + \, \textnormal{HOT}(u) \ \ \ (c_2\ne 0).
\]
\end{lemma}
\noi Writing $C_+, C_- \subseteq A$ for the translates of $C$ by $t_+$ and $t_-$ respectively, the conclusion of the Lemma is illustrated in Figure \ref{AJCurves.Figure}.

 \begin{proof}
 In a neighborhood of $\infty$, $x_\infty$ is expressed as a power series in the local parameter $y_\infty$ with quadratic leading  term:
 \[
 x_\infty \ = \ y_\infty^2 \, + \, \ldots .
 \]
 Therefore $\omega$ and $\eta$ have  local expansions
 \[
 \omega \, = \, (-2 y_\infty^2 + \ldots) \cdot dy_\infty \ \ , \ \ \eta \, = \, (-2 + \ldots) \cdot dy_\infty.
 \] But $\omega = du | C$ and $\eta = dv|C$, so integrating we find that locally
 \[ u \, = \, c \cdot y_\infty^3 + \, \ldots \ \ , \ \ v \, = \, c^\pr \cdot y_\infty \, + \, \ldots
 \]
 with $c, c^\pr \ne 0$. The first statement follows. The assertions for \(t_\pm\) are similar: using \(x\) as local parameter at
\(p_\pm\), one has \(y=\pm1+\ldots\), so
\[
        u-u(t_\pm)=c x+\ldots,\qquad
        v-v(t_\pm)=\pm c' x^2+\ldots
\]
with \(c,c'\ne0\).
 \end{proof}
 
 \parasection{Set-up: The $\mmu_5$-action}  Set $\zeta = \exp(2\pi i/5)$ and consider the action of $\mmu_5 = \langle  \zeta \rangle$ on $C$ defined  by
 \[
 \sigma_C : C \lra C \ \ , \ \ \sigma_C(x,y) \, = \, (\zeta  x, y). 
 \]
 This has the three fixed points $p_\pm , \infty$, and $C / \langle \zeta \rangle = \PP^1$ with the quotient map $C \lra \PP^1$ being totally ramified over three points and \'etale elsewhere. Note that 
 \[ \sigma_C^* (\omega) \ = \  \zeta \cdot \omega \ \ , \ \  \sigma_C^* (\eta) \ = \  \zeta^2 \cdot \eta,\]
 and therefore $\sigma_C$ induces an automorphism
$\sigma : A \lra A$ that in the neighborhood of the origin is given by 
\begin{equation} \label{Local.uv.Action}
\sigma(u, v) \ = \ ( \zeta \cdot u, \zeta^2 \cdot v). 
\end{equation}
This automorphism has five fixed points, namely the multiples of $t_+ \in A$. Set 
\[ t_i \, = \,i \cdot t_+ \ \in \ A \ \   (0 \le i \le 4). \]
 Thus $t_0 = 0$, and in what follows we understand that  the indices are taken (mod $5$). Near each fixed point the action is given in the obvious coordinates by \eqref{Local.uv.Action}. Denote by
 \[ C_i \, = \, C \, + \, t_i  \ \subseteq A\]
 the $t_i$-translate of $C$.\footnote{Thus in Figure 1, $C = C_0, C_+ = C_{1}$ and $C_- = C_{-1}$.} Thus each $C_i$ passes through three of the fixed points, namely $t_{i-1}, t_i$ and $t_{i+1}$. The $C_i$ are stable for the action of $\sigma$, and  $\sigma$ restricts on each of these to $\sigma_C$. 
 
 \parasection{Statement of Main Theorem}
 
 Now consider the quotient map
 \[ \pi : A \lra Y_0 \, =_{\text{def}} \, A/\mmu_5. \]
This is a five-fold cyclic covering, and we fix $\sigma$ as the generator of the Galois group of $\pi$ with which to apply the discussion in \S\ref{Cyclic.Generalities}. Let 
\[
\mu : Y \lra Y_0
\]
be the minimal resolution of the quotient. 
\begin{theorem} \label{Main.Thm.Quotient}
There is a canonically defined birational morphism
\[ Y \lra \PP^1 \times \PP^1. \]
Over $\PP^1 \times \PP^1$, $\pi$ is represented by a cyclic covering branched along a reduced divisor $B \subseteq \PP^1 \times \PP^1$ consisting of the union of five $(1,1)$-curves. 
\end{theorem}
\noi We will see that the birational morphism in question is the composition of two blowings-down of some (-1) curves.

This immediately implies Theorem \ref{Jacobian.Thm.Intro} from the Introduction. In fact, choose a point $o \in \PP^1 \times \PP^1$ away from the branch divisor. Blow up $o$ and then blow down proper transform of the two rulings through $o$ to get a birational map $\PP^1 \times \PP^1 \dra \PP^2$. It follows from Lemma \ref{Branch.Transformation.Rule} that the resulting cover of $\PP^2$ still has reduced branch divisor (consisting of the union of five conics). Hence Proposition \ref{Constructing.Normal.Surfaces} applies. Note however that while the model of the covering over $\PP^1 \times \PP^1$ is canonical, its avatar over $\PP^2$ is not. 
 
\begin{figure} 
\centering
\includegraphics[scale = .7]{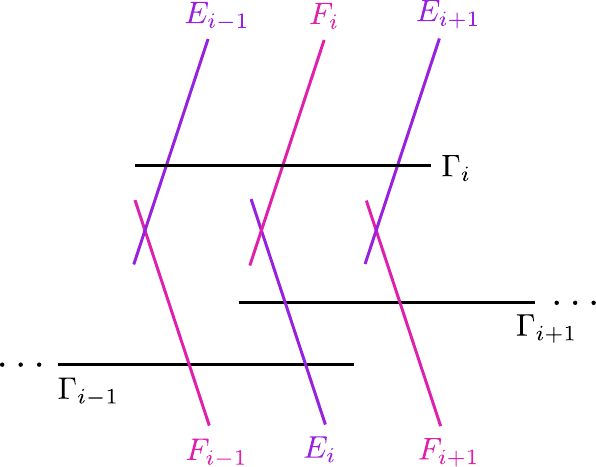}
\caption{Configuration of curves on $Y$}
\label{Curve.Config.on.Y.Figure}
\end{figure}

\parasection{Analysis of the minimal resolution $\mu:  Y \lra Y_0$}\label{Birational.Parasection}
 We next describe the resolution of the quotient, and some curves on it. This is the most crucial part of the analysis.

Each of the fixed points $t_i \in A$ of $\sigma$ gives rise to a singular point  $y_i \in Y_0$, which is local analytically the germ of the quotient of $\CC^2$ by the action \eqref{Local.uv.Action}. The  minimal resolution of these so-called $\tfrac{1}{5}(1,2)$ singularities are well-understood, and Appendix \ref{Appendix} gives a detailed explicit account of the construction.  The upshot is that each of the $t_i$ contributes two exceptional divisors $E_i$, $F_i$ to the resolution $Y$. These are smooth rational curves meeting transversely at a point, with
\begin{equation} \label{EFIntersections}
E_i^2 \, = \, -3 \ \ , \ \ F_i^2 \, = \, -2 .
\end{equation}
By \eqref{Kummer.ef.Multipliciies}, their Kummer branch multiplicities are
\begin{equation} \label{Kummer.Mults.Y}
k_{E_i} \ = 1 \ \ , \ \ k_{F_i}\ = \ 3. 
\end{equation}
Thanks to \eqref{Quotient.Discrepancies} we have 
\begin{equation} \label{Canon.Bundle}
K_Y \ \num \ \ \mu^*(K_{Y_0}) \,- \, \frac{1}{5} \, \sum \, (2E_i + F_i),
\end{equation}
and therefore $5K_Y \lin -\sum (2E_i + F_i)$ thanks to the fact that $K_A = 0$. Since we're aiming to find a model of the covering where all the branch multiplicities $ = 1$, we'd like to contract the $F_i$. It turns out that there is a nice way to do this involving the curves $C_i \subseteq A$. 

Namely, let $\Gamma_i \subseteq Y$ be the proper transform of $\pi(C_i) \subseteq Y_0$. Then 
\[
\Gamma_i \ = \ C_i / \langle 
\sigma \rangle \ = \ \PP^1
\]
is a smooth rational curve, and   \eqref{Canon.Bundle} shows that $\Gamma_i^2 = -1$. Moreover, it follows from Lemma \ref{local.eqn.AJ.Curve} and the computation at the end of the Appendix that $\Gamma_i$ meets $F_i$ transversely at the image of $t_i$, while it meets $E_{i-1}$ and $E_{i+1}$ transversely at the images of $t_{i-1}$ and $t_{i+1}$. Moreover the distinct $\Gamma_i$ are disjoint. This is illustrated schematically in Figure \ref{Curve.Config.on.Y.Figure}. Note that each $F$-curve meets only one $\Gamma$, while each $E$-curve meets two.

\begin{figure} 
\centering
\includegraphics[scale = .73]{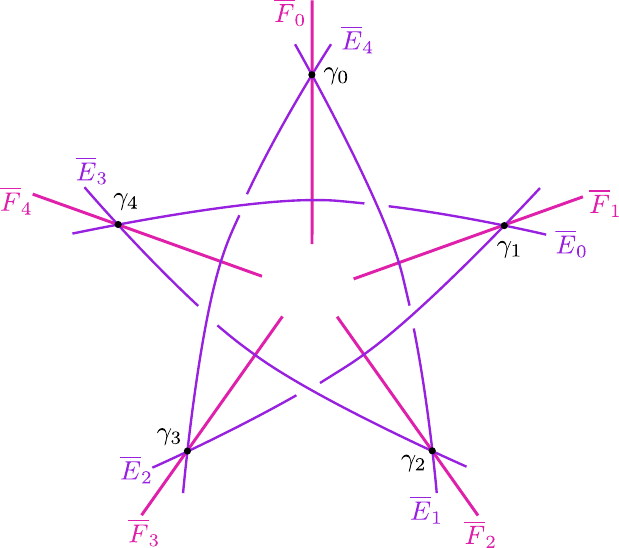}
\caption{Configuration of curves on $\ol{Y}$}
\label{Curve.Config.Figure}
\end{figure}

\parasection{The two contractions}  \label{Two.Contractions} Let
\[
\rho : Y \lra  \ol{Y}
\] be the contraction of the five (-1)-curves $\Gamma_i$. Set $\gamma_i = \rho(\Gamma_i)$  and denote by  $\ol{E}_i, \ol{F}_i \subset  \ol{Y}$ the images of the corresponding curves on $Y$. The configuration of these curves is shown schematically in Figure \ref{Curve.Config.Figure}. On $Y$ each $\Gamma$-curve meets  one $F$-curve but two $E$-curves, which means that on $\ol{Y}$:
\[
\ol{F}_i^2 \, = \, F_i^2 + 1 \, = \, -1 \ \ , \ \  \ol{E}_i^2 \, = \, E_i^2 +2 \, = \, -1.\]
Thanks to Lemma \ref{Branch.Transformation.Rule} the Kummer branch multiplicities of these curves has not changed, ie
\[
k_{\ol{E}_i} \ = 1 \ \ , \ \ k_{\ol{F}_i}\ = \ 3. 
\]

Now blow down the five  $\ol{F}_i$:
\[
\ol{\rho} : \ol{Y} \lra P.
\]
Let
\[ B_i \ = \ \ol{\rho}_*(\ol{E}_i) \ \subseteq \ P. \]
Then $B_i ^2 = 2$ since $\ol{E}_i$ lies on three of the curves being contracted, and $k_{B_i} = 1$.  We see that $B_i$ and $B_{i+2}$ meet a single point with multiplicity two, whereas $B_i$ meets the other $B_j$ transversely at two points.

\parasection{Conclusion} \label{Conclusion.Parasection} We are now reduced to:
\begin{proposition} The surface $P$ is isomorphic to $\PP^1 \times \PP^1$ by an isomorphism under which the $B_i$ are curves of type $(1,1)$.
\end{proposition}

\begin{proof}
Note to begin with that $Y$ is rational thanks to Castelnuovo's criterion. In fact, as $Y_0$ has rational singularities
\[ \HH{1}{Y}{\OO_Y} \ = \ \HH{1}{A}{\OO_A}^{\mmu_5} \ = 0, \]
and since 
$5K_Y \lin -\sum (2E_i + F_i)$ is anti-effective, all its plurigenera vanish.\footnote{Since $P$ is rational and the covering over $P$  has reduced branch divisor, Lemma \ref{Rational.Quotient.Lemma}  already shows that \(A\) admits a normal model in $\PP^3$, but this ignores the interesting geometry of the example.
}
We have $K_Y^2 = -2$, so by Noether's formula $e(Y) = 14$ and hence $\rho(Y) = 12$. Since we obtain $P$ by contracting ten (-1)-curves, this means that $\rho(P) = 2$. Thus \(P\) is a smooth rational surface of Picard number \(2\), and hence
\(P\cong \mathbf F_n\) for some \(n\geq 0\). Note also that $K_P \cdot B_i = -4 $ since $B_i$ is a smooth rational curve of self-intersection = 2.

If $n = 0$ we are done, so assume that $n \ge 1$. Let \(S\) be the negative section on \(P = \mathbb F_n\), and let \(f\) be a
fiber, so that  $K_{\mathbb F_n}\lin -2S-(n+2)f$. Writing
$
        B_i \ \equiv\ aS+bf,
$
and recalling that $K_{\mathbb F_n}\cdot B_i = -4$, we find that 
\[
        B_i \ \lin \  S+\tfrac{n+2}{2}f.
\]
By the standard description of irreducible curves on a Hirzebruch surface, this forces $n\leq 2$. So we must have $n=0$ or $n=2$, and it remains only to rule out the possibility that $n=2$. In fact, if $n=2$ then $B_i\equiv_{\mathrm{lin}}S+2f$, so $B_i\cdot S=0$. Thus $S$ is disjoint from the branch divisor. Moreover, since the images in $P$ of all curves contracted by $Y\to P$ lie on the branch divisor, $S$ lifts to a rational curve in $Y$ disjoint from the exceptional locus of $Y\to Y_0$. Its image in $Y_0$ therefore lies in the smooth locus, and its inverse image in $A$ is a finite étale cover of $\PP^1$, hence a disjoint union of rational curves. This is impossible on an abelian surface.
\end{proof}


\section{A deformation} \label{Deformation.Section}

It turns out that the example analyzed in the previous section deforms in a one-parameter family of genus two Jacobians admitting normal models in $\PP^3$. We give here a brief sketch of how this goes, skipping over many of the details. 

We begin by revisiting the $\mmu_5$ action discussed in \S  \ref{Genus.2.Jacobian.Section}. Referring to Appendix \ref{Appendix}, a local computation shows that the rational quotient mapping
$ \AAA^2_{u,v}  \dra Z^\pr$ can be resolved by blowing up $\AAA^2$ first at the origin and then once in the $v = 0$ direction and twice in the $u = 0$ direction.  It follows that in the global setting the quotient $A \dra Y$ is regularized after the corresponding blowups at each of the fixed points $t_i \in A$. However one finds by tracing through the contractions in Paragraph \ref{Two.Contractions} that the composition
\begin{equation} \label{Quotient.on.A}
A \dra Y \lra \ol{Y} \lra P
\end{equation}
is resolved after just the first two blowups at each of the fixed points.

Starting then with the abelian surface $A$, construct $\nu: A^\pr\lra A$ by first blowing up each of the points $t_i\in A$, creating exceptional divisors $T_i$, and then blowing up each $T_i$ at the point corresponding to the tangent direction $v=0$. Denote by $T_i,R_i\subseteq A^\pr$ the resulting exceptional curves, where $T_i$ now denotes the proper transform of the first exceptional divisor. Thus
\[
T_i^2\, =\, -2\ \ , \ \ R_i^2\, =\, -1.
\]
Write $g : A^\pr \lra P = \PP^1 \times \PP^1$ for the morphism resolving \eqref{Quotient.on.A}.
Next let $D_i\subseteq A^\pr$ denote the proper transform of $C_i$, and set
\begin{equation} \label{Definition.of.H}
H\, = \, \sum_i(D_i+T_i+R_i).
\end{equation}
Tracing the five branch curves $B_i\subseteq P$ backwards through the contractions of the previous section, and using the local description in Appendix \ref{Appendix}, one checks that 
\begin{equation} \label{Pullback.O(11)}
\OO_{A^\pr}(H)\ =\ g^*\OO_{\PP^1 \times \PP^1}(1,1).
\end{equation}
Notice that
\begin{equation} \label{No.Sections.Over.0}
\hh{0}{A^\pr}{H}\ =\ 5:
\end{equation}
we get four sections from $\OO_{\PP^1 \times \PP^1}(1,1)$ and one more from the covering coordinate.

The idea now is that while the $\mmu_5$ action depends on the specific choice of $C$, the configuration of curves on $A$ and $A^\pr$ that we've just discussed -- and ultimately the map to $\PP^1 \times \PP^1$ -- exists on a larger family of Jacobians. Specifically, given $\alpha, \beta \in \CC$, denote by $C_{\alpha, \beta}$ the curve
\[
y^2 \ = \ x^5 \,  + \, (1 + \alpha x + \beta x^2)^2, 
\]
so that $C = C_{0,0}$. For general $\alpha, \beta \in \CC$ this is a smooth hyperelliptic curve, and as before we set $p_{\pm} = (0, \pm 1) \in C_{\alpha, \beta}$. Taking $\omega, \eta$ as in \eqref{Basis.of.Diffs}, Lemma \ref{local.eqn.AJ.Curve} remains valid for the Abel--Jacobi embedding
\[
a \, : \, C_{\alpha, \beta} \hookrightarrow \Jac(C_{\alpha, \beta}), 
\] 
although the specific constants appearing in the local equations depend of course on $\alpha$ and $\beta$.
Moreover, as before $t_+ = a( p_+)$ is a five-torsion point; write $t_i = i \cdot t_+$. While the $\mmu_5$ action is lost, we get a $\ZZ/5\ZZ$ action on $\Jac(C_{\alpha, \beta})$ via translation by $t_+$.

Next, blow up the Jacobian of $C_{\alpha, \beta}$ first at each $t_i$ and then at the points corresponding to the $v =0$ directions to construct
\[
A^\pr_{\alpha, \beta} \lra \Jac(C_{\alpha, \beta}).
\]
Denote by $T_i , R_i \subseteq A^\pr_{\alpha, \beta}$ the corresponding exceptional divisors, and write $D_i \subseteq A^\pr_{\alpha, \beta}$ for the proper transform of the $t_i$-translate of $C_{\alpha, \beta}$.\footnote{We are slightly abusing notation here by using the same letter to denote the corresponding divisors on different $A_{\alpha, \beta}$. We trust that this won't cause any confusion.}  Motivated by the discussion above we set
\[
H\, = \, \sum_i(D_i+T_i+R_i).
\]
The essential point is the following:
\begin{theorem} \label{Section.3.Thm}
There is a Zariski-open neighborhood of the origin in $\CC^2_{\alpha, \beta}$ for which  the line bundle $\OO_{A^\pr_{\alpha, \beta}}(H)$ is basepoint-free and
\[
\hh{0}{A^\pr_{\alpha, \beta}}{\OO_{A^\pr_{\alpha, \beta}}(H)} \ = \ 5.
\]
The   resulting morphism 
\[
A^\pr_{\alpha, \beta} \lra \PP^4
\]
 is birational onto its image, which  lies on a quadric. Moreover, there is a curve 
 \[
 (0, 0) \, \in \, \Delta \, \subseteq \,  \CC^2_{\alpha, \beta}\]
 such that for general \((\alpha,\beta)\in\Delta\) near $(0,0)$, the corresponding quadric has rank \(4\), and the image of $A'_{\alpha,\beta}$ is a  subvariety of the total space of $\mathcal O_{\mathbf P^1\times\mathbf P^1}(1,1)$.
 \end{theorem}

The Theorem implies that if  $(\alpha, \beta) \in \Delta$ is a general point, then $\Jac(C_{\alpha, \beta})$ is birational to a normal decic $S_{\alpha, \beta} \subseteq \PP^3$. In fact, let 
\[ W_{\alpha, \beta} \ \subseteq \ \text{Tot}\big( \OO_{\PP^1 \times \PP^1}(1,1) \big)\]
denote the image of $A^\pr_{\alpha, \beta}$. Thus $W_{0,0}$ is the cyclic cover analyzed in the previous section. It follows that $W_{\alpha,  \beta}$ is a monic hypersurface of relative degree $5$, and these images vary in a flat family of divisors. Since $W_{0,0}$ is normal, so also is $W_{\alpha, \beta}$ for  general $(\alpha, \beta) \in \Delta$. Now as in \S \ref{Genus.2.Jacobian.Section}, projection from a general point of $\PP^1 \times \PP^1$ gives rise to a birational isomorphism 
\[ \Tot \big(   \OO_{\PP^1 \times \PP^1}(1,1) \big) \ \sim_{\text{bir}} \ \Tot \big( \OO_{\PP^2}(2) \big). \]
As long as the center of the projection is chosen away from the branch divisor of $W_{\alpha, \beta} \lra \PP^1 \times \PP^1$, its image
\[ V_{\alpha, \beta} \ \subseteq \ \Tot \big( \OO_{\PP^2}(2) \big) \]
is a normal monic hypersurface of relative degree 5.\footnote{For the normality, observe that away from the two points corresponding to the contracted rulings, the transformed covering is either isomorphic to the original covering or étale along the exceptional line.}  Thus we are in the situation of Remark \ref{Monic.Hypersurface} and we get the asserted realization of $\Jac(C_{\alpha, \beta})$.

\begin{proof}[Sketch of Proof of Theorem \ref{Section.3.Thm}]
We first show that $H$ has at least five independent sections for every smooth curve $C_{\alpha,\beta}$
in the family.  

The local calculations of Lemma
\ref{local.eqn.AJ.Curve} show that the intersection configuration of these fifteen curves is unchanged
under our deformation. In particular,
\[
D_i^2\, =\, -3\ \ , \ \  T_i^2\, = \, -2\ \ , \ \ R_i^2\, = \, -1.
\]
Observe that the $D_i$ are disjoint. Furthermore, $T_i$ meets precisely $D_i$ and $R_i$, while $R_i$ meets $T_i,D_{i-1},D_{i+1}$ at three
distinct points. 

We claim first that
\[
h^0\big(H,\OO_H(H)\big)\ =\ 5.
\]
In fact, adjunction and the intersection configuration give
\[
\OO_{D_i}(H)
\, \cong \, K_{C_{\alpha,\beta}}(-p_+-p_-)
\, \cong\, \OO_{C_{\alpha,\beta}},
\]
the last equality holding because $p_++p_-$ is a
canonical divisor. Similarly,
\[
\OO_{T_i}(H)\, \cong \,\OO_{\PP^1},
\qquad
\OO_{R_i}(H)\, \cong\, \OO_{\PP^1}(2).
\]
We can now extend an arbitrary choice of sections
\[
s_i \, \in \, \HH{0}{D_i}{\OO_{D_i}(H)}\ \cong\ \CC\]
to all of $H$. The $s_i$ agree with unique sections of $\OO_{T_i}(H)$, and hence prescribe values at the three attaching points of each $R_i$. Since a section of $\OO_{\PP^1}(2)$ is specified by its values at three distinct points, we then obtain a section of $\OO_H(H)$. Conversely every such section arises in this way, proving the claim.  

Next consider the exact sequence
\[
0\lra\OO_{A^\pr_{\alpha,\beta}}
\lra\OO_{A^\pr_{\alpha,\beta}}(H)
\lra\OO_H(H)\lra0.
\]
Translation by $t_+$ lifts to $A^\pr_{\alpha,\beta}$
and cyclically permutes the five triples
$(D_i,T_i,R_i)$. The discussion in the previous paragraph shows that the induced action of $\ZZ/5\ZZ$ on $H^0(H,\OO_H(H))$ is the regular
representation. On the other hand, translation acts trivially on
\[
\HH{1}{A^\pr_{\alpha,\beta}}{\OO_{A^\pr_{\alpha,\beta}}}
\ =\ 
\HH{1}{\Jac(C_{\alpha,\beta})}{
       \OO_{\Jac(C_{\alpha,\beta})}}
\]
The connecting homomorphism
\[
H^0(H,\OO_H(H))
\lra H^1(A^\pr_{\alpha,\beta},
              \OO_{A^\pr_{\alpha,\beta}})
\]
is equivariant, and therefore factors through
the one-dimensional invariant subspace of
$H^0(H,\OO_H(H))$. Its rank is consequently
at most one, and we conclude that
\[
h^0(A^\pr_{\alpha,\beta},
             \OO_{A^\pr_{\alpha,\beta}}(H)) \tag{*}
\, \geq \, 1+5-1\, =\, 5.
\]

But now return to the central fibre $A^\pr = A^\pr_{0,0}$. As noted in \eqref{Pullback.O(11)} and \eqref{No.Sections.Over.0}, there we know that
$
\hh{0}{A^\pr}{H} = 5
$
and that  $\OO_{A^\pr}(H)=g^*\OO_P(1,1)$.
In the first place, it follows by semicontinuity from (*) that
\[
\hh{0}{A^\pr_{\alpha, \beta}}{H} \,  = \, 5
\]
for general $\alpha, \beta$. Moreover, the number of sections being constant, the basepoint-freeness of 
$\OO_{A^\pr}(H)$ implies the same for $(\alpha, \beta)$ in a neighborhood of $(0,0)$. 
 We consequently obtain morphisms
\[
\varphi_{\alpha,\beta}:
A^\pr_{\alpha,\beta}\lra\PP^4
\]
for general $(\alpha, \beta)$.

We next show that the images  of the $\varphi_{\alpha, \beta}$ lie on a quadric. On the central fibre, the cyclic covering of
$P=\PP^1\times\PP^1$ gives
\[
h^0(A^\pr,2H)=14.
\]
Moreover, the multiplication map
\[
\Sym^2 H^0(A^\pr,H)
\lra H^0(A^\pr,2H)
\]
is surjective,  its one-dimensional kernel being generated by the equation of the rank-four
quadric containing the central image. Fourteen independent quadratic products remain
independent under deformation, while by semicontinuity $\hh{0}{A^\pr_{\alpha,\beta}}{2H} \le 14$. It follows
that the multiplication map continues to have rank fourteen, and hence that every nearby image
lies on a unique quadric $Q_{\alpha,\beta}=0$.

Finally, observe that $Q_{0,0}$ has rank four: its zero-locus is the (one-point compactification of) the total space of $\OO_{\PP^1 \times \PP^1}(1,1)$. Thus there is a curve $\Delta \subseteq \CC^2_{\alpha, \beta}$ containing $(0,0)$ parameterizing points for which $\rk( Q_{\alpha, \beta}) = 4$.   Moreover, for general $(\alpha, \beta) \in \Delta$ the image of $\varphi_{\alpha, \beta}$ misses the vertex of $Q_{\alpha, \beta}$, since this is true for $(0,0)$, and since $\varphi_{0,0}$ is birational to its image so too is the general $\varphi_{\alpha, \beta}$. \end{proof}

\begin{remark}
It is perhaps worth noting  that the family $\{ C_{\alpha, \beta}\}$ is non-isotrivial. In fact, a fixed genus two curve has only finitely many choices of a Weierstrass point $\infty$ and a point $p_+$ satisfying $5[p_+-\infty]=0$, and once these are fixed, the normalization of the equation of $C_{\alpha,\beta}$ leaves only finitely many possibilities for $(\alpha,\beta)$.
\end{remark}

\appendix  

\section{Resolution of the $\frac15(1,2)$ cyclic quotient singularity} \label{Appendix}

We give here an explicit description, in coordinates, of the resolution
of the cyclic quotient singularity that plays a central role above.
These calculations can be used to verify the geometric assertions
appearing in the body of the paper.

Consider the action of $\mmu_5$ on $\mathbf A^2_{u,v}$ given by
\begin{equation} \label{Action.on.C2}
 \zeta\cdot(u,v)\ =\ (\zeta u,\zeta^2v),
 \qquad \zeta \, \in \, \mmu_5 \ \text{ a generator}.
\end{equation}
Let
\[
Z\ =  \ \mathbf{A}^2/\mmu_5
\]
be the quotient. Its coordinate ring is the algebra generated by the
invariant polynomials
\[
 Z\ =\ \operatorname{Spec}\bigl(\CC[u^5\, , \, u^3v\, , \, uv^2\, , \, v^5]\bigr).
\]
We  fix the automorphism \eqref{Action.on.C2} as a generator of the Galois group of the quotient.
The minimal resolution\[
 \mu\colon Z'\longrightarrow Z
\]
can be constructed by two weighted blow-ups that create   exceptional curves
\begin{equation} \label{EF.Facts}
E \ , \ F \ \subseteq \ Z',
 \qquad E^2\, = \, -3,\quad F^2\, = \, -2,
\end{equation}
both isomorphic to $\PP^1$, that meet  transversely at a point. But it is more convenient to work concretely. Specifically,   $Z^\pr$ is covered by three copies 
of $\mathbf A^2$:
\[
 U_L\, =\, \mathbf A^2_{a,b},
 \qquad U_M\, = \, \mathbf A^2_{e,f},
 \qquad U_R\, = \, \mathbf A^2_{c,d},
\]
determined by the inclusions shown in Figure \ref{fig:quotient-resolution-charts}. These three charts glue via the evident identifications to give the surface $Z^\pr$  and the map $Z^\pr \lra Z$.

\begin{figure}[t]
\[
\xymatrix@C=1.4em@R=-.45em{
 \CC[a,b] && \CC[e,f] && \CC[c,d] \\
 \raisebox{0pt}[\dimexpr\height+.95em\relax][\dimexpr\depth+.35em\relax]{\rotatebox{90}{\scalebox{1.25}[1]{$=$}}} &&
 \raisebox{0pt}[\dimexpr\height+.95em\relax][\dimexpr\depth+.35em\relax]{\rotatebox{90}{\scalebox{1.25}[1]{$=$}}} &&
 \raisebox{0pt}[\dimexpr\height+.95em\relax][\dimexpr\depth+.35em\relax]{\rotatebox{90}{\scalebox{1.25}[1]{$=$}}} \\
 \CC\left[u^5,\dfrac{v}{u^2}\right] &&
 \CC\left[\dfrac{v^3}{u},\dfrac{u^2}{v}\right] &&
 \CC\left[\dfrac{u}{v^3},v^5\right] \\
 & \raisebox{.45em}{\rotatebox{-35}{\large$\supseteq$}} &
 \rule[-1.6em]{0pt}{3.2em}\raisebox{.45em}{\rotatebox{-90}{\large$\supseteq$}} &
 \raisebox{.45em}{\rotatebox{35}{\large$\subseteq$}} & \\
 && \CC[u^5,u^3v,uv^2,v^5] &&
}
\]
\caption{The three affine charts in the resolution of the singularity
\(\frac{1}{5}(1,2)\).}
\label{fig:quotient-resolution-charts}
\end{figure}

The divisors $E$ and $F$ are defined on these three charts by the local equations:
\[
\renewcommand{\arraystretch}{1.35}
\begin{array}{c|ccc}
   & U_L & U_M & U_R \\ \cline{1-4}
E  & a=0 & e=0 & {}  \\[4pt]
F  & {}  & f=0 & d=0
\end{array}
\]  The transitions are given by
\[
e\, = \, ab^3,\qquad f\, = \, b^{-1},\qquad c\, = \, e^{-1},\qquad d\, = \, e^2f,
\]
and the assertions in \eqref{EF.Facts} can be deduced from this. For example, along $E$  the elements $b$ and $f=b^{-1}$ are affine coordinates on the two standard charts of $\PP^1$, while the normal coordinates satisfy $e=b^3a$. Hence
\[
N_{E/Z'}\cong \mathcal O_{\PP^1}(-3),
\]
and similarly the relation $d=e^2f$ gives
\[
N_{F/Z'}\cong \mathcal O_{\PP^1}(-2).
\]
Applying the adjunction formula to  \eqref{EF.Facts}, a standard computation of discrepancies shows that 
\begin{equation} \label{Quotient.Discrepancies}
    \ord_E(K_{Z^\pr} )\ = \ - \tfrac{2}{5} \ \ , \ \  \ord_F(K_{Z^\pr}) \ = \ - \tfrac{1}{5}.
\end{equation}

The invariant functions in each chart are given by the table:
\begin{equation} \label{Local.Generator.Table}
\renewcommand{\arraystretch}{1.35}
\begin{array}{c|cccc}
 & u^5 & u^3v & uv^2 & v^5 \\ \hline
U_L
 & a
 & ab
 & ab^2
 & a^2b^5
\\[4pt]
U_M
 & ef^3
 & ef^2
 & ef
 & e^2f
\\[4pt]
U_R
 & c^5d^3
 & c^3d^2
 & cd
 & d
\end{array}
\end{equation}
Having fixed \eqref{Action.on.C2} as a generator of the Galois group of the quotient, the function $u^5$ is the Kummer generator of the branch divisor in the sense of \S \ref{Cyclic.Generalities}. Hence we read off from \eqref{Local.Generator.Table} the Kummer multiplicities of the two exceptional divisors:
\begin{equation} \label{Kummer.ef.Multipliciies}
k_E \ = \ 1 \ \ , \ \ k_F \, = \, 3. 
\end{equation}

Finally, suppose that  $\gamma\subseteq\CC^2_{u,v}$ is an
invariant curve with local analytic equation at the origin
\[
 u\ =\ c_1   v^3\, + \, \operatorname{HOT}(v),
 \qquad c_1 \, \ne \, 0. \]
 Then in the first place, $\mu_5$-invariance forces the higher-order terms to have order $\ge 8$. Moreover, consider the proper transform $\gamma^\pr \subseteq Z^\pr$ of the image of $\gamma$ under the quotient map $\CC_{u,v}^2 
 \lra Z$. Then the local equations show that 
 $\gamma^\pr$ in $Z^\pr$ meets $F$
transversely. Similarly, the proper transform of the image of an invariant curve
\[
 v\ = \ c_2 u^2\, + \, \operatorname{HOT}(u) 
\]
meets $E$ transversely, and different values of $c_2$ lead to different
intersection points.


\end{document}